\documentclass{hjm1}
\usepackage{verbatim}
\usepackage{amssymb}

   \theoremstyle{remark}
   \newtheorem{rem}[thm]{Remark}
  \theoremstyle{definition}
   \newtheorem{example}[thm]{Example}

\newcommand{\topq}{E=(E^{0},E^{1},r,s,\lambda)}
\newcommand{\supp}{\operatorname{supp}}

\begin{document}

\title{Markov Operators and $C^{*}$-Algebras }

\author{Marius Ionescu}

\address{Department of Mathematics, 196 Auditorium Road, Unit 3009 University
of Connecticut, Storrs, CT 06269-3009}

\email{ionescu@math.uconn.edu}

\author{Paul S. Muhly}

\address{Department of Mathematics, The University of Iowa, Iowa City, IA
52242-1419}

\email{pmuhly@math.uiowa.edu}

\author{Victor Vega}

\address{Department of Mathematics, St. Ambrose University, 518 W. Locust
St, Davenport, IA 52803}

\email{VegaVictorM@sau.edu}
\begin{abstract}
A Markov operator $P$ acting on $C(X)$, where $X$ is compact, gives
rise to a natural topological quiver. We use the theory of such quivers
to attach a $C^{*}$-algebra to $P$ in a fashion that reflects some
of the probabilistic properties of $P$. 
\end{abstract}
\maketitle

\section{Introduction}

\global\long\def\supp{\operatorname{supp}}
Our objective in this note is to use the theory of topological quivers
\cite{MuSo_PLMS00,MuTo05} to study natural $C^{*}$-algebras that
can be associated to Markov operators. In particular, we shall use
the theory developed in \cite{MuTo05} to decide when these $C^{*}$-algebras
are simple. In addition, we will explore a number of examples that
help illustrate how our analysis may be applied and we shall explore
connections between Markov operators and topological quivers. 

The term ``Markov operator'' is used in a variety senses in the probability
literature. We adopt the following definition here and will discuss
the terminology more in Remark \ref{rem:Transpose P}.
\begin{dfn}
\label{def:Markov-operator}Let $X$ be a compact Hausdorff space
and let $C(X)$ denote the space of continuous, complex-valued functions
on $X$. A \emph{Markov operator} on $C(X)$ is a unital positive
linear map $P$ on $C(X)$.
\end{dfn}
%
{}
\begin{dfn}
\label{def:topological-quiver}A \emph{topological quiver} is a quintet
$\topq$, where $E^{0}$ and $E^{1}$ are second countable, locally
compact Hausdorff spaces, $r$ and $s$ are continuous maps from $E^{1}$
onto $E^{0}$, with $r$ open, and where $\lambda=\{\lambda_{v}\}_{v\in E^{0}}$
is a family of measures on $E^{1}$such that the (closed) support
of $\lambda_{v}$, $\supp\lambda_{v}$, equals $r^{-1}(v)$, and such
that for each function $f\in C_{c}(E^{1})$, the function $v\to\int_{E^{1}}f(x)\, d\lambda_{v}$
lies in $C_{c}(E^{0})$. The space $E^{0}$ is called the space of
\emph{vertices} of $E$, $E^{1}$ is the space of \emph{edges}, $r$
and $s$ are called the \emph{range} and \emph{source} maps, respectively,
and $\lambda$ is called the family of weights. 
\end{dfn}
In a fashion that will be spelled out in a bit more detail in the
next section, each Markov operator $P$ on $C(X)$ gives rise to a
topological quiver $E$. The vertex space $E^{0}$ is $X$, the edge
space $E^{1}$ is the ``support'' of $P$, a subspace of $X\times X$,
the range and source maps are the left and right projections, respectively,
and the family of weights is given by a continuous family of probability
measures naturally associated to $P$. The $C^{*}$-algebra that we
associate to $P$ and will denote by $\mathcal{O}(P)$ is the $C^{*}$-algebra
of this quiver.

In the next section, we provide additional definitions and detail,
and we provide a variety of examples (not exhaustive) to which our
analysis applies. Section \ref{sec:Simplicity} is devoted to determining
when $\mathcal{O}(P)$ is simple. Section \ref{sec:Examples} is devoted
to applying our simplicity criteria to the examples described in Section
\ref{sec:Definitions-and-Examples}. Finally, in Section \ref{sec:From-Quivers}
we address the question: Given a topological quiver $\topq$ when
can one find a Markov operator $P$ so that the $C^{*}$-algebra of
$\topq$ is isomorphic to $\mathcal{O}(P)$?

\section{Definitions and Examples\label{sec:Definitions-and-Examples}}

Throughout this note, $X$ will be a fixed compact Hausdorff space,
which we shall assume to be second countable. Also, $P$ will be a
fixed Markov operator on $C(X)$. The Riesz representation theorem
gives a continuous family of measures on $X$ and indexed by $X$,
which we will write $p(\cdot,y)$, such that \begin{equation}
(Pf)(y)=\int f(x)p(dx,y).\label{eq:measures}\end{equation}
The fact that $P$ is unital implies that each $p(\cdot,y)$ is a
probability measure. Further, we can always extend $P$ to the bounded
Borel functions on $X$ via the formula $Pf(y):=\int f(x)p(dx,y)$
and, consequently, $p(U,y)=\int1_{U}(x)p(dx,y)=P(1_{U})(y)$ for all
Borel sets $U\subseteq X$ and all $y\in X$. The natural topological
quiver to associate to $P$ is the one that gives the so-called GNS
correspondence for $P$. Our first objective in this section is to
give details to support this assertion.
\begin{rem}
\label{rem:Transpose P}Before continuing we want to explain why some
may find our formula for $P$ unconventional. At one point early in
the theory, Markov operators were defined as certain operators acting
on measures on measurable spaces. Various hypotheses were imposed
to insure that the operators had ``adjoints'' that acted on the space
of measurable functions. The formula for the adjoint action on functions
was written {\begin{equation}
(Pf)(x)=\int p(x,dy)f(y).\label{eq:Transposea}\end{equation}
Evidently, \eqref{eq:Transposea} is a transposed version of \eqref{eq:measures}.
We have chosen our notation, \eqref{eq:measures}, to conform with
the conventions from graph algebra theory and other situations where
one builds operator algebras from not-necessarily-reversable dynamical
systems. As a result, some of our formulas are transposed versions
of formulas in the literature. We note, too, that in the probability
literature, what we are calling a Markov operator is sometimes called
a Markov-\emph{Feller} operator. Feller identified conditions that
insure that $P$, defined on measurable functions as in \eqref{eq:measures}
or \eqref{eq:Transposea}, leaves the space of continuous functions
invariant (assuming, of course, the measure space is built on some
topological space.) For a contemporary view of these issues, and references
see \cite{rZ2005}.}
\end{rem}
It is well-known that since $C(X)$ is commutative, $P$ is completely
positive. Consequently, the following definition makes sense. 
\begin{dfn}
\label{def:GNS}The \emph{$GNS$-correspondence} for $P$ (over the
$C^{*}$-algebra $C(X)$), is the space $C(X)\otimes_{P}C(X)$, which
is the separated completion of the algebraic tensor product $C(X)\odot C(X)$
in the pre-inner product defined by the formula $\langle\xi_{1}\otimes\eta_{1},\xi_{2}\otimes\eta_{2}\rangle:=\overline{\eta_{1}}P(\overline{\xi_{1}}\xi_{2})\eta_{2}$,
and is endowed with the bimodule structure over $C(X)$ defined by
the formula $a\cdot(\xi\otimes\eta)\cdot b:=(a\xi)\otimes(\eta b)$. 
\end{dfn}
{}
\begin{dfn}
\label{def:support}The \emph{support} of $P$ (or of $p$), denoted
$\supp(P)$ (resp. $\supp(p)$) is the complement of the set of all
points $(x,y)\in X\times X$ with the property that there is a neighborhood
$U$ of $x$ such that the function $y\to p(U,y)$ vanishes in some
neighborhood of $y$. Equivalently, $(x_{0},y_{0})\notin\supp P$
if and only if there is a neighborhood of $(x_{0},y_{0})$ of the
form $U\times V$ such that $p(U,y)=0$, for all $y\in V$. 
\end{dfn}
%
{}
\begin{prop}
\label{pro:MO-to-TopQuiv}Let $E^{0}$ be $X$, let $E^{1}\subseteq X\times X$
be $\supp(P)$, define $r$ and $s$ by the formulae $s(x,y)=x$ and
$r(x,y)=y$, respectively, and define $\lambda_{v}:=p(\cdot,v)$,
$v\in E^{0}$. Then the quintet $\topq$ is a topological quiver,
with $E^{1}$compact. Further, $C(E^{1})$ becomes a $C^{*}$-correspondence
$\mathcal{X}$ over $C(X)$ via the formula\[
a\cdot\xi\cdot b(x,y)=a(x)\xi(x,y)b(y),\]
 and\[
\langle\xi,\eta\rangle(y)=\int\overline{\xi(x,y)}\eta(x,y)p(dx,y)=\int\overline{\xi}\eta\, d\lambda_{y},\]
 $\xi,\eta\in C(E^{1})$, $a,b\in C(E^{0})=C(X)$, and the map $W:C(X)\otimes_{P}C(X)\to\mathcal{X}$
defined by $W(f\otimes g)(x,y)=f(x)g(y)$ extends to an isomorphism
of correspondences from the $GNS$-correspondence for $P$ to $\mathcal{X}$.\end{prop}
\begin{proof}
With all the definitions before us, the proof is nothing but a straightforward
process of checking. Of course $E^{1}$ is compact, since it is a
closed subset of $X\times X$. The other matters are equally easy. 
\end{proof}
The correspondence $\mathcal{X}$ is the correspondence of the kind
that is associated to any topological quiver \cite[Subsection 3.1]{MuTo05}.
\begin{dfn}
\label{def:E(P)-X(P)} The topological quiver associated to $P$ in
Proposition \ref{pro:MO-to-TopQuiv} will be denoted $E(P)$ and the
resulting correspondence will be denoted $\mathcal{X}(P)$.\end{dfn}
\begin{rem}
The topological quiver $E(P)$ is a topological relation in the sense
of Brenken \cite{BrB04} since $E^{1}$ is a closed subset of $X\times X$.
\end{rem}
To define the $C^{*}$-algebra that we associate to $P$, and to relate
it to $E(P)$, it will be helpful to spell out additional definitions
and facts that will also play a role elsewhere in this note. Given
a $C^{*}$-correspondence $\mathcal{X}$ over a $C^{*}$-algebra $A$
a \emph{Toeplitz representation} of $\mathcal{X}$ in a $C^{*}$-algebra
$B$ consists of a pair $(\psi,\pi)$, where $\psi:\mathcal{X}\rightarrow B$
is a linear map and $\pi:A\rightarrow B$ is a $\ast$-homomorphism
such that\[
\psi(x\cdot a)=\psi(x)\pi(a)\;,\;\psi(a\cdot x)=\pi(a)\psi(x),\]
 i.e. the pair $(\psi,\pi)$ is a \emph{bimodule map}, and such that\[
\psi(x)^{*}\psi(y)=\pi(\langle x,y\rangle_{A}).\]
 That is, the map $\psi$ preserves inner products (see \emph{\cite[Section 1]{FoMuRa03}}).
Given such a Toeplitz representation, there is a $\ast$-homomorphism
$\pi^{(1)}$ from $\mathcal{K}(\mathcal{X})$ into $B$ which satisfies\begin{equation}
\pi^{(1)}(\Theta_{x,y})=\psi(x)\psi(y)^{*}\;\;\mbox{for}\;\mbox{all}\; x,y\in\mathcal{X},\label{eq:PiSuper1}\end{equation}
 where $\Theta_{x,y}=x\otimes\tilde{y}$ is the rank one operator
defined by $\Theta_{x,y}(z)=x\cdot\langle y,z\rangle_{A}$. The $*$-homomorphism
extends naturally to $\mathcal{L}(\mathcal{X})$ by virtue of the
fact that $\mathcal{L}(\mathcal{X})$ is the multiplier algebra of
$\mathcal{K}(\mathcal{X})$ and we will denote the extension by $\pi^{(1)}$
also. This extension $\pi^{(1)}$ and $\psi$ are related by the useful
formula \begin{equation}
\pi^{(1)}(T)\psi(\xi)=\psi(T\xi),\label{eq:Pi1ExtAction}\end{equation}
$T\in\mathcal{L}(\mathcal{X})$ and $\xi\in\mathcal{X}$. Indeed,
if $T$ is a rank one operator, $\Theta_{x,y}$, then \[
\pi^{(1)}(T)\psi(\xi)=\psi(x)\psi(y)^{*}\psi(\xi)=\psi(x)\pi(\langle y,\xi\rangle)=\psi(\Theta_{x,y}(\xi)).\]
So the formula holds for all $T\in\mathcal{K}(\mathcal{X})$. If $T\in\mathcal{L}(\mathcal{X})$
and $S\in\mathcal{K}(\mathcal{X})$, then\begin{eqnarray*}
\pi^{(1)}(T)\psi(S\xi) & = & \pi^{(1)}(T)\pi^{(1)}(S)\psi(\xi)=\pi^{(1)}(TS)\psi(\xi)\\
 & = & \psi((TS)\xi)=\psi(T(S\xi)),\end{eqnarray*}
since $\mathcal{K}(\mathcal{X})$ is an ideal in $\mathcal{L}(\mathcal{X})$
and $\mathcal{L}(\mathcal{X})$ is the multiplier algebra of $\mathcal{K}(\mathcal{X})$.
Thus $\pi^{(1)}(T)\psi(\xi)=\psi(T\xi)$ for all $T\in\mathcal{L}(\mathcal{X})$
and $\xi\in\mathcal{X}$.

Let $\Phi:A\to\mathcal{L}(\mathcal{X})$ be the $*$-homomorphism
that defines the left action of $A$ on $\mathcal{X}$. We define
then \[
J(\mathcal{X}):=\Phi^{-1}(\mathcal{K}(\mathcal{X})),\]
 which is a closed two sided-ideal in $A$ (see \cite[Definition 1.1]{FoMuRa03})
and we set $J_{\mathcal{X}}:=J(\mathcal{X})\bigcap\ker^{\perp}\Phi$,
where $\ker^{\perp}\Phi$ denotes the set of all $a\in A$ such that
$ab=0$ for all $b\in\ker\Phi$. Suppose $K$ is any ideal in $J(\mathcal{X})$.
We say that a Toeplitz representation $(\psi,\pi)$ of $\mathcal{X}$
is \emph{coisometric} \emph{on $K$} if\[
\pi^{(1)}(\Phi(a))=\pi(a)\;\mbox{for}\;\mbox{all}\; a\in K.\]
 When $(\psi,\pi)$ is coisometric on all of $J(\mathcal{X})$, we
say that it is \emph{Cuntz-Pimsner covariant.}

It is shown in \cite[Proposition 1.3]{FoMuRa03} that for an ideal
$K$ in $J(\mathcal{X})$, there is a $C^{*}$-algebra $\mathcal{O}(K,\mathcal{X})$,
called the \emph{relative Cuntz Pimsner algebra associated to $\mathcal{X}$
and $K$}, and a Toeplitz representation $(k_{\mathcal{X}},k_{A})$
of $\mathcal{X}$ into $\mathcal{O}(K,\mathcal{X})$, which is coisometric
on \emph{$K$,} and satisfies:
\begin{enumerate}
\item for every Toeplitz representation $(\psi,\pi)$ of $\mathcal{X}$
which is coisometric on $K$, there is a $\ast$-homomorphism $\psi\times_{K}\pi$
of $\mathcal{O}(K,\mathcal{X})$ such that $(\psi\times_{K}\pi)\circ k_{\mathcal{X}}=\psi$
and $(\psi\times_{K}\pi)\circ k_{A}=\pi$; and
\item $\mathcal{O}(K,\mathcal{X})$ is generated as a $C^{*}$-algebra by
$k_{\mathcal{X}}(\mathcal{X})\cup k_{A}(A)$. 
\end{enumerate}
The algebra $\mathcal{O}(\{0\},\mathcal{X})$, then, is the \emph{Toeplitz
algebra} $\mathcal{T}_{\mathcal{X}}$, and $\mathcal{O}(J_{\mathcal{X}},\mathcal{X})$
is defined to be the \emph{Cuntz-Pimsner algebra $\mathcal{O}_{\mathcal{X}}$
}(see \cite{tKa03} and \cite{MuTo05})\emph{.}

The parts of the following definition are taken from \cite[Section 3]{MuTo05}
\begin{dfn}
\label{def:Vertices-Revue}Let $E=(E^{0},E^{1},r,s,\lambda)$ be a
topological quiver and let $\mathcal{X}$ be the $C^{*}$-correspondence
over $A=C(E^{0})$ associated to $E$. \end{dfn}
\begin{enumerate}
\item The set of \emph{sinks} of $E$, $E_{sinks}^{0}$, is defined to be
the open subset $U$ of $E^{0}$ that supports the kernel of $\Phi$,
i.e., $U$ satisfies the equation $\Phi^{-1}(0)=C_{0}(U)$. 
\item The set of \emph{finite emitters} of $E$, $E_{fin}^{0}$, is defined
to be the open subset of $E^{0}$ that supports $\Phi^{-1}(\mathcal{K}(\mathcal{X}))$,
i.e., $\Phi^{-1}(\mathcal{K}(\mathcal{X}))=C_{0}(E_{fin}^{0})$. 
\item A vertex $v$ is called regular if it is a finite emitter, but not
a sink. The set of all regular vertices is denoted $E_{reg}$, so
that $E_{reg}^{0}=E_{fin}^{0}\setminus E_{sinks}^{0}$. 
\item Elements of $E^{0}\setminus E_{fin}^{0}$are called \emph{infinite
emitters. }
\item The $C^{*}$-algebra $C^{*}(E)$ associated to $E$ is defined to
be the relative Cuntz-Pimsner algebra $\mathcal{O}(C_{0}(E_{reg}^{0}),\mathcal{X})$.
\end{enumerate}
{}

{}
\begin{rem}
Since the Markov operator $P$ is unital, $E(P)$ has no \emph{sources},
that is $r(E^{1})=E^{0}$. Otherwise, if $x\in E^{0}\setminus r(E^{1})$
then $P(1)(x)=0$, by a compactness argument, which is a contradiction.
Further, since $s(E^{1})$ is compact, and hence closed, $E_{sinks}=E^{0}\setminus s(E^{1})$.
Moreover $x$ is a sink if and only if there is an open neighborhood
$U$ of $x$ such that $P(1_{U})(y)=0$ for all $y\in E^{0}$. Indeed,
if $x$ is sink then for each $y\in E^{0}$ there is a neighborhood
$U_{y}$ of $x$ and a neighborhood $V_{y}$ of $y$ such that $z\mapsto p(U_{y},z)$
is $0$ on $V_{y}$. Since $E^{0}$ is compact, there is a finite
subcover $\{V_{1},V_{2},\dots,V_{n}\}$ of $X$. Let $U=\bigcap_{i=1}^{n}U_{i}$,
which is open. It follows that $z\mapsto p(U,z)=0$ for all $z\in E^{0}$,
that is, $P(1_{U})=0$. The converse is clear. Finally, note that
because $P$ is unital, $\mathcal{X}(P)$ is \emph{full}, meaning
that the ideal in $C(X)$ generated by the inner products equals $C(X)$.\end{rem}
\begin{dfn}
The $C^{*}$-algebra of  $P$ is defined to be
$\mbox{ }C^{*}(E(P))$ and will be denoted $\mathcal{O}(P)$.\end{dfn}
\begin{example}
\label{exa:FiniteMC}If $X=\{x_{1},\dots,x_{n}\}$ is a finite set,
then $P$ is a Markov operator on $X$ if and only if there is a stochastic
matrix $p$ such that\[
P(f)(x_{j})=\sum_{i}f(x_{i})p_{ij},\]
for all $f\in C(X)$. To say $p$ is stochastic means here that $\sum_{i}p_{ij}=1$,
which is the transpose of the usual definition (see, for example \cite{KeSn_FMC60}).
In this case, of course, $E(P)$ is a finite directed graph with vertices
$E^{0}=X$, and the set of edges consists of the pairs $(x_{i},x_{j})$
such that $p_{ij}>0$. Consequently, $\mathcal{O}(P)$ is the graph
$C^{*}$-algebra studied by the third author in \cite{Vega_PHD}.
More accurately, he initially defined $\mathcal{O}(P)$ to be the
Cuntz-Pimsner algebra of the $GNS$-correspondence determined by $P$
and arrived at the representation of $\mathcal{O}(P)$ as a graph
$C^{*}$-algebra through (a special case of) Proposition \ref{pro:MO-to-TopQuiv}.
\end{example}
{}
\begin{example}
\label{exa:homeom}If $\tau:X\to X$ is a homeomorphism, and if $p(\cdot,y)$
is the point mass at $\tau^{-1}(y)$, then $P$ is the automorphism
$\alpha$ of $C(X)$ given by the formula $\alpha(f)(y)=f(\tau^{-1}(y))$
since \[
Pf(y)=\int f(x)p(dx,y)=\int\delta_{\tau^{-1}(y)}(dx)f(x)=f(\tau^{-1}(y))=\alpha(f)(y).\]
The support of $P$ is the graph of $\tau$. Thus the topological
quiver is $E(P)=(E^{0},E^{1},r,s,\lambda)$ where $E^{0}=X$, $E^{1}$
is the graph of $\tau$, i.e. $E^{1}=\{(x,\tau(x))\,:\, x\in X\}$,
$r(x,\tau(x))=\tau(x)$, $s(x,\tau(x))=x$, and\begin{eqnarray*}
a\cdot\xi\cdot b(x,\tau(x)) & = & a(x)\xi(x,\tau(x))b(\tau(x)),\\
\langle\xi,\eta\rangle(y) & = & \overline{\xi(\tau^{-1}(y),y)}\eta(\tau^{-1}(y),y).\end{eqnarray*}
The $C^{*}$-algebra $\mathcal{O}(P)$ is then the cross-product $C^{*}$-algebra
$C(X)\rtimes_{\alpha}\mathbb{Z}$. Of course, viewing crossed products
as Cuntz-Pimsner algebras was one of Pimsner's sources of inspiration
\cite{mP97}.
\end{example}
{}
\begin{example}
\label{exa:loc_hom} More generally, let $\tau:X\to X$ be a local
homeomorphism. Then $\tau^{-1}(y)$ is a finite set for all $y\in X$.
If $p(\cdot,y)$ is counting measure on $\tau^{-1}(y)$ normalized
to have total mass $1$, then\[
P(f)(y)=\frac{1}{\vert\tau^{-1}(y)\vert}\sum_{x\in\tau^{-1}(y)}f(x).\]
The support of $P$ is still the graph of $\tau$. The only difference
from the previous example is the formula for the inner product, which
becomes\[
\langle\xi,\eta\rangle(y)=\frac{1}{\vert\tau^{-1}(y)\vert}\sum_{x\in\tau^{-1}(y)}\overline{\xi(x,y)}\eta(x,y).\]
The $C^{*}$-algebra $\mathcal{O}(P)$ is the cross-product of $C(X)$
by the local homeomorphism $\tau$ studied extensively in \cite{Exel_ETDS03,BroRae_MPCPS06,Exel_JFA03,ExVe_CJM06,ArzRen_96,Deac_TAMS95,Deac_PAMS96,Deac_PJM99,Deac_Theta00,DeKuMu,KuRe_PAMS06,IoWa_CMB08,ExLo_ETDS04}.
\end{example}
{}
\begin{example}
\label{exa:ifs}Markov operators have played a role in the theory
of iterated function systems from the very begining of the theory
fractals. See Hutchinson's paper \cite{Hutc_80} where they are mentioned
explicitly in this context. As Hutchinson notes, many of the ideas
he develops can be traced back further in geometric measure theory.
Barnsley and his collaborators used Markov operators to good effect
in encoding and decoding pictures in terms of fractals. For a sampling
of this literature, see \cite{mB93,mB2006,BE88}. If $X$ is a compact
metric space, an \emph{iterated function system} on $X$ is a finite
set of injective contractions $(f_{1},f_{2},\dots,f_{N})$, $N\ge2$,
on $X$. Given such an iterated function system there is a unique
compact subset $K$ of $X$ which is invariant for the iterated function
system, that is\[
K=f_{1}(K)\bigcup f_{2}(K)\bigcup\dots\bigcup f_{n}(K).\]
In the following we assume that $X=K$. Then\[
P(f)(y):=\frac{1}{N}\sum_{i=1}^{N}f\circ f_{i}(y)\]
is a Markov operator on $C(X)$.

\global\long\def\cograph{\operatorname{cograph}}
We claim that the support of $P$ equals $\bigcup_{i=1}^{N}\cograph f_{i}$,
where for a function $f:X\to X$ the \emph{cograph} of $f$ is\[
\cograph f=\{(x,y)\,:\, x=f(y)\}.\]
Indeed, if $(x,y)$ belongs to the support of $P$ then for any neighborhoods
$U$ and $V$ of $x$ and $y$, respectively, there exists $z\in V$
and $i\in\{1,\dots,N\}$ such that $f_{i}(z)\in U$. Therefore we
can find a sequence $\{z_{n}\}_{n}$ that converges to $y$ and a
sequence of indices $\{i_{n}\}_{n}\subset\{1,\dots,N\}$ such that
$\lim_{n\to\infty}f_{i_{n}}(z_{n})=x$. There must be an index $i\in\{1,\dots,N\}$
so that $i_{n}=i$ infinitely many times. Therefore there is a subsequence
$\{z_{n_{k}}\}$ so that $\lim_{k\to\infty}f_{i}(z_{n_{k}})=x$. Thus
$x=f_{i}(y)$ and $(x,y)$ belongs to the cograph of $f_{i}$. The
converse inclusion is clear. Thus the topological quiver $E(P)$ is
given by $E^{0}=X$, $E^{1}=\bigcup_{i=1}^{N}\cograph f_{i}$, $s(x,y)=x$,
$r(x,y)=y$. The actions and inner product on $\mathcal{X}(P)$ are
given by\begin{eqnarray*}
(a\cdot\xi\cdot b)(x,y) & = & a(x)\xi(x,y)b(y)\\
\langle\xi,\eta\rangle(y) & = & \frac{1}{N}\sum_{i=1}^{N}\overline{\xi(f_{i}(y),y)}\eta(f_{i}(y),y).\end{eqnarray*}
Then the $C^{*}$-algebra $\mathcal{O}(P)$ is the $C^{*}$-algebra
studied in \cite{KaWa_JOT06} and in \cite{IoWa_CMB08}.

More generally, if $p=\{p_{1},p_{2},\dots,p_{N}\}$ are probabilities
($p_{i}>0$ for all $i$ and $\sum_{i=1}^{N}p_{i}=1$), then $P_{p}(f)(y):=\sum_{i=1}^{N}p_{i}f\circ f_{i}(y)$
is a Markov operator on $C(X)$. Since $p_{i}>0$ for all $i\in\{1,\dots,N\}$
it follows that the support of $P_{p}$ is still $\bigcup_{i=1}^{N}\cograph f_{i}$.
The only difference between $\mathcal{X}(P)$ and the $C^{*}$-correspondence
$\mathcal{X}(P_{p})$ associated to $P_{p}$ is the formula for the
inner product, which is\[
\langle\xi,\eta\rangle_{p}(y)=\sum_{i=1}^{N}p_{i}\overline{\xi(f_{i}(y),y)}\eta(f_{i}(y),y).\]
We claim that $\mathcal{X}(P)$ and $\mathcal{X}(P_{p})$ are isomorphic
$C^{*}$-correspondences (\cite{MuSo_PLMS00}). To prove the claim
recall first that for $(x,y)$ in the support of $P$ (which is the
same as the support of $P_{p}$) we define it's branch index to be
$e(x,y)=\#\{i\in\{1,\dots,N\}\,:\, f_{i}(y)=x\}$ (see \cite{KaWa_JOT06}
and \cite{IoWa_CMB08}). Then one can prove that the map $\psi:\mathcal{X}(P)\to\mathcal{X}(P_{p})$
defined by\[
\psi(\xi)(x,y)=\frac{e(x,y)^{1/2}}{\sqrt{N}\bigl(\sum_{i\,:\, f_{i}(y)=x}p_{i}\bigr)^{1/2}}\xi(x,y),\]
for $\xi\in C(E^{1})$, is a $C^{*}$-correspondence isomorphism.
Thus $\mathcal{O}(P)$ and $\mathcal{O}(P_{p})$ are isomorphic as
$C^{*}$-algebras. For simplicity we will assume in the sequel that
$p_{1}=p_{2}=\dots=p_{N}=1/N$.
\end{example}
{}
\begin{example}
\label{exa:nonIFS}Let $X=[0,1]$, let $f_{1}(x)=x$, and let $f_{2}(x)=1-x$.
Then $f_{1}$ and $f_{2}$ are not contractions, so $(f_{1},f_{2})$
is not an iterated function system. Nevertheless,\[
P(f)(y):=\frac{1}{2}\sum_{i=1}^{2}f\circ f_{i}(y)\]
is a Markov operator. More generally, if $(f_{1},f_{2},\dots,f_{N})$
are continuous functions on $X$ so that $X=f_{1}(X)\bigcup f_{2}(X)\bigcup\dots\bigcup f_{N}(X)$,
then \[
P(f)(y):=\frac{1}{N}\sum_{i=1}^{N}f\circ f_{i}(y)\]
is a Markov operator on $C(X)$. The support of $P$ is still $\bigcup_{i=1}^{N}\cograph f_{i}$.
\end{example}
{}
\begin{example}
\label{exa:non_rowfinite}If $X=[0,1]$ and $p(\cdot,y)=m$, Lebesgue
measure, for all $y\in X$, then $P$ is given by the formula\[
P(f)(y):=\int_{[0,1]}f(x)dm(x)\]
for all $y\in X$. The support of $P$ is $X\times X$, $E^{0}=X$,
$E^{1}=X\times X$, $r(x,y)=y$, $s(x,y)=x$, $\lambda_{y}=m$ for
all $y\in X$. The actions and inner product on the associated $C^{*}$-correspondence
are given by the formulae \begin{eqnarray*}
a\cdot\xi\cdot b(x,y) & = & a(x)\xi(x,y)b(y),\\
\langle\xi,\eta\rangle(y) & = & \int\overline{\xi(x,y)}\eta(x,y)dm(x).\end{eqnarray*}
This is the prototypical example of a topological quiver for which
every vertex is an infinite emitter.
\end{example}

\section{Simplicity of the $C^{*}$-algebras\label{sec:Simplicity}}

Given a topological quiver $E=(E^{0},E^{1},r,s,\lambda)$, the $C^{*}$-algebra
$C^{*}(E)$ is simple if and only if $E$ satisfies condition (L)
and the only open saturated hereditary subsets of $E^{0}$ are $E^{0}$
and $\emptyset$ \cite[Theorem 10.2]{MuTo05}. In this section we
investigate these conditions for the topological quivers derived from
Markov operators.%
{} For this we need to review a few definitions (see, for example \cite{MuTo05,tKa04,tKa06_2}).

A \emph{path of length $n$} in a topological quiver $E=(E^{0},E^{1},r,s,\lambda)$
is a finite sequence $\alpha=\alpha_{1}\alpha_{2}\dots\alpha_{n}$
so that $\alpha_{i}\in E^{1}$ and $r(\alpha_{i})=s(\alpha_{i+1})$
for all $i=1,\dots,n-1$. We denote by $E^{n}$ the set of paths of
length $n$, we write $E^{*}:=\bigcup_{n\ge0}E^{n}$ for the set of
all \emph{finite} paths, and we write \[
E^{\infty}:=\{\alpha=(\alpha_{n})_{n\ge0}\,:\,\alpha_{n}\in E^{1}\mbox{ and }r(\alpha_{n})=s(\alpha_{n+1})\}\]
for the set of \emph{infinite paths}. For a finite path $\alpha=\alpha_{1}\alpha_{2}\dots\alpha_{n}$
we define $s(\alpha):=s(\alpha_{1})$ and $r(\alpha):=r(\alpha_{n})$.
If $\alpha$ is an infinite path we define $s(\alpha):=s(\alpha_{1})$.
A finite path $\alpha$ is called a \emph{loop }if $r(\alpha)=s(\alpha)$.
In this case we say that $v=r(\alpha)=s(\alpha)$ is a \emph{base}
of the loop $\alpha$. If $\alpha=\alpha_{1}\alpha_{2}\dots\alpha_{n}$
is a loop, then an \emph{exit for $\alpha$} is an edge $\beta\in E^{1}$
such that $s(\beta)=s(\alpha_{i})$ for some $i\in\{1,\dots,n\}$
and $\beta\ne\alpha_{i}$ (\cite[Definition 6.8]{MuTo05}). Thus,
if $x=s(\alpha_{i})$, there must be at least two edges leaving $x$,
i.e., $\vert s^{-1}(x)\vert\ge2$. 
\begin{dfn}
\label{def:Cond-L-Top_Quiv} A topological quiver is said to satisfy\emph{
condition (L)} (in the sense of \cite[Definition 6.9]{MuTo05}) if
the set of base points of loops with no exit has empty interior.
\end{dfn}
%
{}Turning now to condition (L) for our Markov operator $P$, observe
that $P^{n}$ is also Markov operator on $C(X)$ for all $n\ge1$.
In the next lemma we relate the paths of $E$ to the powers of $P$.
%
{}
\begin{prop}
\label{pro:path_power}For $v,w\in X$ there is a path of length $n$
from $v$ to $w$ if and only if for any neighborhood $V$ of $v$
and neighborhood $W$ of $w$ there is $f\in C(X)$ with $\supp f\subset V$
and $f\ge0$, such that $P^{n}(f)(z)>0$ for some $z\in W$.\end{prop}
\begin{proof}
In one direction, we proceed by induction based on the length of the
path from $v$ to $w$. If the length of the path is one, the condition
means that $(v,w)$ belongs to the support of $P$. The conclusion
is, then, the definition of the support of $P$. Assume next that
there is a path of length two from $v$ to $w$, that is, assume there
is an $\alpha=\alpha_{1}\alpha_{2}$, $\alpha_{i}\in E^{1}$, such
that $s(\alpha)=v$, $r(\alpha)=w$. Let $V$ be an open neighborhood
of $v$ ad $W$ an open neighborhood of $w$. Then for any open neighborhood
$V_{1}$ of $r(\alpha_{1})=s(\alpha_{2})$ there are $f\in C(X)$
with $\supp f\subset V$ and $f\ge0$, and $f_{1}\in C(X)$ with $\supp f_{1}\subset V_{1}$
and $0\le f_{1}\le1$, such that $P(f_{1}P(f))(z)>0$ for some $z\in W$.
Then\begin{eqnarray*}
P^{2}(f)(z) & = & \int_{X}P(f)(y)p(dy,z)\ge\int_{X}f_{1}(y)P(f)(y)p(dy,z)\\
 & = & P(f_{1}P(f))(z)>0.\end{eqnarray*}
The inductive step is now clear. The other direction is immediate.
\end{proof}
Next we consider saturated hereditary sets for Markov operators basing
the terminology on \cite[Definition 8.3]{MuTo05}. 
\begin{dfn}
\label{def:Hereditary-MO}If $(E^{0},E^{1},r,s,\lambda)$ is topological
quiver, then a subset $U\subseteq E^{0}$ is \emph{hereditary }if
whenever $\alpha\in E^{1}$ and $s(\alpha)\in U$, then $r(\alpha)\in U$.
A subset $U$ of $X$ will be called hereditary for the Markov operator
$P$ in case $U$ is hereditary for $E(P)$. An hereditary subset
$U$ of $E^{0}$ is called \emph{saturated }if whenever $x\in E_{reg}^{0}$
and $r(s^{-1}(x))\subseteq U$, then $x\in U$.\end{dfn}
\begin{lem}
\label{lem:support}Suppose $y\in X$ and there is an open set $U$
such that $\int_{U}p(dx,y)>0$. Then there is an $x\in U$ such that
$(x,y)$ belongs to the support of $P$.\end{lem}
\begin{proof}
Since $X$ is assumed to be a second countable compact space, $X$
is metrizable. Therefore, we may choose a complete metric on $X$
so that the topology is given by the metric. %
{} Since $p(\cdot,y)$ is a Radon measure there is $x_{1}\in U$ and
$r_{1}>0$ so that if $V_{1}=B(x_{1},r_{1})$ - the ball of radius
$r_{1}$, centered at $x_{1}$, computed with respect to the metric
- then $\int_{V_{1}}p(dx,y)>0$ and $V_{1}\subseteq\overline{V}_{1}\subseteq U$.
Inductively we may find a sequence of points $\{x_{n}\}$ and radii
$\{r_{n}\}$ such that $r_{n}\to0$ and if $V_{n}:=B(x_{n},r_{n})$,
then $V_{n}\supset V_{n+1}$ and $\int_{V_{n}}p(dx,y)>0$ for all
$n\ge1.$ The sequence $\{x_{n}\}_{n\in\mathbb{N}}$ is Cauchy and
so we may let $x:=\lim_{n\to\infty}x_{n}$. Then from construction
it follows that $x\in U$ and $\int_{V}p(dx,y)>0$ for all neighborhoods
$V$ of $x$. Thus $(x,y)$ belongs to the support of $P$.
\end{proof}
{}
\begin{prop}
\label{pro:Markov_hereditary}If $U$ is an open subset of $X$ that
is hereditary for the Markov operator $P$, then $P$ restricts to
a Markov operator $P_{K}$ on $C(K)$, where $K=X\setminus U$.\end{prop}
\begin{proof}
Recall that the quotient $C(X)/C_{0}(U)$ is $*$-isomorphic with
$C(K)$ via the map $[f]\to f|_{K}$, where $[f]$ is the equivalence
class of an element $f\in C(X)$ in $C(X)/C_{0}(U)$, and $f|_{K}$
is the restriction of $f$ to $K$. Therefore it is enough to check
that if $U$ is hereditary then $P$ maps $C_{0}(U)$ into $C_{0}(U)$.
This is clear, however, since if $\int_{U}p(dx,y)>0$ for some $y\in K$
it follows that that there is $x\in U$ such that $(x,y)\in E^{1}$
by Lemma \ref{lem:support}. Then $y\in U$ since $U$ is assumed
to be hereditary, and this is a contradiction.
\end{proof}
%
{}

Recall from \cite[Definition 8.8]{MuTo05} that if $\mathcal{E}=(E^{0},E^{1},r,s,\lambda)$
is a topological quiver and $U$ is a hereditary open subset of $E^{0}$,
then one can ``cut $\mathcal{E}$ down to $U$'' to obtain a topological
quiver $\mathcal{E}_{U}:=(E_{U}^{0},E_{U}^{1},r_{U},s_{U},\lambda_{U})$,
where $E_{U}^{0}=E^{0}\setminus U$, $E_{U}^{1}=E^{1}\setminus r^{-1}(U)$,
$r_{U}=r|_{E_{U}^{1}}$, $s_{U}=s|_{E_{U}^{1}}$, and $\lambda_{U}=\lambda|_{E_{U}^{0}}$.
Thus, if $\mathcal{E}$ is the topological quiver associated with
a Markov operator and $U$ is a hereditary set, then $\mathcal{E}_{U}$
is the topological quiver associated with $P_{K}$, where $K=X\setminus U$
and $P_{K}$ is the restriction of $P$ to $C(K)$. Recall also the
following definition from the theory of Markov chains (see, e.g.,
\cite[Definition 2.2.2]{HerLass_MC}).
\begin{dfn}
\label{def:Absorbing} A Borel subset $B$ of $X$ is called \emph{absorbing}
with respect to the Markov operator $P$, if $p(B,x)=1$ whenever
$x$ lies in $B$. 
\end{dfn}
Thus, $B$ is absorbing for $P$ if and only if $P(1_{B})$ is identically
$1$ on $B$.
\begin{prop}
\label{pro:Ssaturated}If $U$ is an open hereditary set, then $K:=X\setminus U$
is a closed absorbing set.\end{prop}
\begin{proof}
%
{}If $K$ is not an absorbing set, then there is a $y\in K$ such that
$\int_{U}p(dx,y)>0$. Lemma \ref{lem:support} implies that there
is an $x\in U$ such that $(x,y)$ belongs to the support of $P$.
Since $U$ is hereditary it follows that $y\in U$, which is a contradiction.\end{proof}
\begin{prop}
A closed set $K$ is absorbing set for $P$ if and only if $P$ restricts
to a Markov operator $P_{K}$ on C(K).\end{prop}
\begin{proof}
Suppose $K$ is an absorbing set for $P$. Then $P(1_{K})(x)=1$ for
all $x\in K$. It follows that if $f\in C_{0}(U)$ then $P(f)(x)=0$
for all $x\in K$. Thus $P(f)\in C_{0}(U)$. Hence $P$ restricts
to a Markov operator on $C(K)$. Conversely, if $P$ restricts to
a Markov operator $P_{K}$ on $C(K)$, then for any $f\in C(X)$,
with $f|_{K}=1$, we have $P_{K}([f])(x)=1$ for all $x\in K$. Thus
$p(K,x)=1$ for all $x\in K$, and $K$ is absorbing for $P$.\end{proof}
\begin{lem}
\label{lem:absorbing_power}If $B$ is an absorbing set for $P$,
then $B$ is absorbing for $P^{n}$ for all $n\ge1$.\end{lem}
\begin{proof}
The proof is by induction. By hypothesis, $P(1_{B})(z)=1$ for all
$z\in B$. Therefore $P(1_{B})-1_{B}\ge0$. Since $P$ is positive,
\[
1\ge P^{2}(1_{B})(z)=P(P(1_{B}))(z)\ge P(1_{B})(z)=1\]
for all $z\in B$. Thus $B$ is absorbing for $P^{2}$. The induction
is now clear.\end{proof}
\begin{dfn}
We say that a Borel set $B$ is \emph{strongly absorbing} for the
Markov operator $P$ if the following condition holds: A point $x$
belongs to the complement of $B$ if and only if there is an open
neighborhood $V$ of $x$ and a an open neighborhood $W$ of $B$
so that $\int_{V}p(dx,y)=0$ for all $y\in W$.
\end{dfn}
Note that a closed strongly absorbing set for $P$ is absorbing. Indeed,
if we fix $y\in B$, we can cover the complement of $B$ with a countable
collection of open sets $\{V_{n}\}_{n\ge0}$ so that $\int_{V_{n}}p(dx,y)=0$
for all $n\in\mathbb{N}$. Then $\int_{B^{c}}p(dx,y)=0$ so $\int_{B}p(dx,y)=1$.
Thus $B$ is absorbing.

We come now to the main result of this section.
\begin{thm}
\label{thm:Simplicity}Suppose %
{}$P$ is a Markov operator on $C(X)$ such that $E(P)$ satisfies condition
(L). %
{} Then $\mathcal{O}(P)$ is simple if and only if the only closed strongly
absorbing set is $X$.\end{thm}
\begin{proof}
According to Theorem 10.2 of \cite{MuTo05} to show that $\mathcal{O}(P)$
is simple, we must show that the only open saturated hereditary sets
in $X$ are $X$ and $\emptyset$. Consequently, it is enough to show
that a nonempty open subset $U$ of $X$ is a saturated hereditary
subset of $E^{0}$ if and only if $K:=X\setminus U$ is a closed strongly
absorbing set for $P$. 

Assume first that $U$ is an open saturated hereditary set. Then since
$K$ is closed, $K$ is compact. Further, since $U$ is hereditary,
$r(s^{-1}(U))\subset U$. So, if $x\in U$, then for each $y\in K$
there is an open neighborhood $V_{y}$ of $x$ and an open neighborhood
$W_{y}$ of $y$ so that $\int_{V_{y}}p(dx,z)=0$ for all $z\in W_{y}$.
Consequently, $\{W_{y}\}_{y\in K}$ is an open cover of $K$. Let
$\{W_{y_{1}},W_{y_{2}},\dots,W_{y_{n}}\}$ be a finite subcover of
$K$ and let $\{V_{y_{1}},V_{y_{2}},\dots,V_{y_{n}}\}$ be the corresponding
open neighborhoods of $x$. Then for $V=\bigcap_{i=1}^{n}V_{y_{i}}$
and $W=\bigcup_{i=1}^{n}W_{y_{i}}$ ,%
{}\[
\int_{V}p(dx,y)=0\;\mbox{for all}\; y\in W.\]
%
{} Thus $K$ is strongly absorbing. %
{}

Assume, conversely, that $K$ is a strongly absorbing subset of $X$.
If $U$ is not hereditary, there is $x\in U$ and $y\in K$ so that
$(x,y)\in\supp P$. Therefore for any neighborhood $V_{x}$ of $x$
and any neighborhood $W_{y}$ of $y$ there is $z\in W_{y}$ so that
$\int_{V_{x}}p(dx,y)>0$. This contradicts the definition of a strongly
absorbing set. Thus $U$ is hereditary. To see that $U$ is saturated,
let $x$ be a regular vertex with the property that $r(s^{-1}(x))\subset U$.
Then for any $y\in K$ there is a neighborhood $V_{y}$ of $x$ and
a neighborhood $W_{y}$ of $y$ so that $\int_{V_{y}}p(dx,z)=0$ for
all $z\in V_{y}$. Since $E_{reg}^{0}=E_{fin}^{0}-E_{sinks}$, we
may assume without loss of generality that each of the $V_{y}$ is
contained in $E_{reg}^{0}$. By the compactness of $K$ again, we
can find a finite number of points such that $\{W_{y_{1}},W_{y_{2}},\dots,W_{y_{n}}\}$
is a cover of $K$. If $V=\bigcap_{i=1}^{n}V_{y_{i}}$ as before,
we obtain an open neighborhood $V$ of $x$ and an open neighborhood
$W=\bigcup_{i=1}^{n}W_{y_{i}}$ of $K$ so that $\int_{V}p(dx,y)=0$
for all $y\in W$. Since $K$ is strongly absorbing it follows that
$x\in U$, thus $U$ is an open saturated hereditary set and the theorem
is proved.
\end{proof}

\section{Examples\label{sec:Examples}}

\subsection{Finite Markov Chains}

Let $X$ be a finite set and $P$ be a Markov operator on $C(X)$.
Recall from Example \ref{exa:FiniteMC} that $P$ is given by (the
transpose of) a stochastic matrix $\{p_{ij}\}$. Then $\mathcal{O}(P)$
is simple if and only if the Markov chain is ergodic (\cite{KeSn_FMC60}).
We recover, thus, Theorem 5.16 of \cite{Vega_PHD}.

{}

\subsection{Homeomorphisms and local homeomorphism}

Let $X$ be a compact second countable Hausdorff space and let $\tau$
be a homeomorphism or a local homeomorphism on $X$. The canonical
Markov operators we defined in Examples \ref{exa:homeom} and \ref{exa:loc_hom}
are defined by \[
P(f)(y)=f(\tau^{-1}(y))\]
 if $\tau$ is a homeomorphism, and \[
P(f)(y)=\frac{1}{\vert\tau^{-1}(y)\vert}\sum_{x\in\tau^{-1}(y)}f(x)\]
if $\tau$ is a local homeomorphism. In both cases the support of
$P$ is the graph of $\tau$. Let $E(P)=(E^{0},E^{1},r,s,\lambda)$
be the associated topological quiver. We describe next a characterization
of condition (L) for this class of examples. We begin with a straightforward
application of Proposition \ref{pro:path_power}.
\begin{cor}
If $P$ is the Markov operator associated with a homeomorphism or
a local homeomorphism $\tau$ on $X$, and $x,y\in X$, then there
is a path from $x$ to $y$ if and only if $y=\tau^{n}(x)$ for some
$n\ge1$. Thus $x$ is a base point of a loop if and only if $x$
is a fixed point for $\tau^{n}$ for some $n\ge1$.
\end{cor}
Since $(x,y)\in E^{1}$ if and only if $\tau(x)=y$, it follows that
for any point $x\in E^{0}$ there is exactly one edge whose source
is $x$. In particular no loop has an exit. Thus we obtain the following
characterization for condition (L).
\begin{cor}
Let $\tau$ be a homeomorphism or a local homeomorphism on $X$. Let
$P$ be the associated Markov operator. Then the topological quiver
$E(P)$ satisfies condition (L) if and only if the set of fixed points
of $\tau^{n}$, $n\ge1$, has empty interior.
\end{cor}
One can easily see that a subset $U$ of $X$ is hereditary for $E$
if and only if $\tau^{-1}(K)\supset K$, where $K=X\setminus U$.
Moreover $U$ is a saturated hereditary subset of $X$ if and only
if $K=\tau^{-1}(K)$. Thus a set $K$ is strongly absorbing if and
only if $K$ is invariant under $\tau$, that is $K=\tau^{-1}(K)$.
Then we obtain the following characterization for the simplicity of
these $C^{*}$-algebras, which in the case of homeomorphisms, at least,
has been known for quite some time.
\begin{prop}
If $P$ is the Markov operator associated to a homeomorphism or local
homeomorphism $\tau$ on $X$, then $\mathcal{O}(P)$ is simple if
and only if the set of fixed points of all positive powers of $\tau$
has empty interior and there are no closed invariant sets under $\tau$.
\end{prop}

\subsection{Iterated Function Systems}

Let $X$ be a compact metric space. Recall from Example \ref{exa:ifs}
that an iterated function system (i.f.s.) is a collection of injective
contractions $\{f_{1},f_{2},\dots,f_{N}\}$ with $N\ge2$. We assume
that $X$ is invariant for the i.f.s., that is\[
X=f_{1}(X)\bigcup f_{2}(X)\bigcup\dots\bigcup f_{N}(X).\]
The canonical Markov operator associated with an i.f.s. is \[
P(f)(y)=\frac{1}{N}\sum_{i=1}^{N}f\circ f_{i}(y)\]
and its support is $\bigcup_{i=1}^{N}\cograph f_{i}$. The associated
topological quiver is given by $E^{0}=X$, $E^{1}=\bigcup_{i=1}^{N}\cograph f_{i}$,
$s(x,y)=x$, $r(x,y)=y$. We will show that these topological quivers
will always satisfy condition (L) (Proposition \ref{pro:IFS_CondL}).
We start with an immediate consequence of Proposition \ref{pro:path_power}.
\begin{cor}
\label{cor:basepoints_ifs}If $P$ is the Markov operator associated
with an iterated function system $(f_{1},f_{2},\dots,f_{N})$ on $X$,
and $x,y\in X$ then there is a path from $x$ to $y$ if and only
if there is a finite word $w=w_{1}w_{2}\dots w_{n}$, $w_{i}\in\{1,\dots,N\}$,
such that $y=f_{w}(x)$, where $f_{\omega}=f_{\omega_{1}}\circ\dots\circ f_{\omega_{n}}$.
Thus $x$ is a base point for a loop if and only if $x$ is the fixed
point of $f_{\omega}$, for some finite word $\omega\in\{1,\dots,N\}^{n}$.\end{cor}
\begin{proof}
For $f\in C(X)$ we have that\[
P^{n}(f)(y)=\frac{1}{N^{n}}\sum_{w\in\{1,\dots,N\}^{n}}f\circ f_{w}(y).\]
Now Proposition \ref{pro:path_power} implies the conclusion.
\end{proof}
Since $X$ is the invariant set of the iterated function system it
follows that if $x\in X$ then there is some $y\in X$ and $i\in\{1,\dots,N\}$
such that $x=f_{i}(y)$. Therefore $(x,y)$ belongs to the support
of $P$ and there is at least one edge with source $x$. We identify
next the classes of iterated function systems for which there is exactly
one edge with source $x$ for all $x\in X$. These i.f.s. are classified
as follows. 

The iterated function system is called \emph{totally disconnected}
if $f_{i}(X)\bigcap f_{j}(X)=\emptyset$ for $i\ne j$. In this case,
there is a local homeomorphism $\tau:X\to X$ such that $\tau\circ f_{i}=1_{X}$
and the Markov operator $P$ is the same as the Markov operator associated
to the local homeomorphism $\tau$.

A point $x\in X$ is called a \emph{branch point (see \cite[Definition 2.4]{KaWa_JOT06},
\cite[Definition 2.4]{IoWa_CMB08})} if there are two indices $i\ne j$
and $y\in X$ such that $x=f_{i}(y)=f_{j}(y)$. If $f_{i}(X)\bigcap f_{j}(X)$
is not empty but consists of only branch points, then there still
is a continuous map $\tau$ such that $\tau\circ f_{i}=1_{K}$. The
map $\tau$ is a branch covering in this case.
\begin{prop}
Let $(f_{1},f_{2},\dots,f_{n})$ be an iterated function system with
invariant set $X$ and let $P$ the associated Markov operator on
$C(X)$. Then for any point $x\in E^{0}$ there is exactly one edge
whose source is $x$ if and only if the iterated function system is
totally disconnected or if $f_{i}(X)\bigcap f_{j}(X)\ne\emptyset$
then it contains only branch points. \end{prop}
\begin{proof}
If the iterated function system is totally disconnected, given any
$x\in X$ there is a unique $y\in X$ and $i\in\{1,\dots,N\}$ such
that $x=f_{i}(y)$. That is, there is a unique $y\in X$ with $(x,y)\in\mbox{supp}P$.
If the iterated function system is not totally disconnected and $f_{i}(X)\bigcap f_{j}(X)$
is either empty or contains branch points, then for any $x\in X$
there is a \emph{unique} $y\in X$ such that $x=f_{i}(y)$ for some
index $i$ which might not be unique. Then $(x,y)$ is the unique
edge in $E^{1}$ with source $x$.

Conversely, if for any $x\in X$ there is a unique edge in $E^{1}$
with source $x$, then there is a unique $y\in X$ such that $(x,y)$
belongs to the support of $P$. Therefore for each $x$ there is a
unique $y$ so that $x=f_{i}(y)$ for some $i\in\{1,\dots,N\}$. This
clearly implies that the iterated function system is either totally
disconnected or $f_{i}(X)\bigcap f_{j}(X)$ consists only of branch
points, if it is nonempty.\end{proof}
\begin{example}
Let $X=[0,1]$, $f_{1}(x)=\frac{1}{2}x$, $f_{2}(x)=1-\frac{1}{2}x$.
Then the union of the cographs of $f_{i}$ is the graph of the tent
map (\cite[Example 4.5]{KaWa_JOT06}) and $f_{1}(X)\bigcap f_{2}(X)$
contains only $1/2$, which is a branch point. Thus for each point
$x\in[0,1]$ there is a unique edge with source $x$. Note, however,
that $1/2$ is \emph{not} a finite emitter in the sense of \cite[Definition 3.14]{MuTo05}.
According to \cite[Proposition 2.6]{KaWa_JOT06}, \cite[Proposition 2.6]{IoWa_CMB08},
the set of finite emitters for this example is $[0,1]\setminus\{1/2\}$.
\end{example}
It follows that if the iterated function system is either totally disconnected
or  $f_{i}(X)\bigcap f_{j}(X)$ contain only branch points ($i\ne j$), no loop
has an exit. Nevertheless, the associated topological quiver will
always satisfy condition (L).
\begin{prop}
\label{pro:IFS_CondL}Assume that $X$ is a non-discrete uncountable
compact metric space and $(f_{1},f_{2},\dots,f_{N})$ is an iterated
function system on $X$. Let \begin{equation}
Pf(y)=\frac{1}{N}\sum_{i}f\circ f_{i}(y)\label{eq:MarkovofIFS}\end{equation}
be the associated Markov operator. Then the corresponding topological
quiver $E(P)$ satisfies condition (L).\end{prop}
\begin{proof}
We will prove that the set of base point of loops has empty interior.
Recall from Lemma \ref{cor:basepoints_ifs} that $x$ is the base
point of a loop in $X$ if and only if $x$ is a fixed point of $f_{w}$
for a finite word $w\in\{1,\dots,N\}^{n}$, for some $n\ge1$. Since
each $f_{w}$ is a contraction it has a unique fixed point. Since
the number of finite words over $\{1,\dots,N\}$ is countable so is
the set of base points of loops, hence it has empty interior.
\end{proof}
Next we study the open saturated hereditary subsets of $E^{0}$ for
an iterated function system. The lack of the existence of such sets
is an example of the {}``rigidity'' of iterated function systems
in operator algebras. For other results suggesting this rigidity see
\cite{PiWaYo_CMP00} and \cite{mIo_RMJM07}.
\begin{prop}
\label{pro:IFS_hereditary}Let $(f_{1},f_{2},\dots,f_{N})$ be an
iterated function system with invariant set $X$. Then the only two
open hereditary subsets are $X$ and the empty set. \end{prop}
\begin{proof}
Suppose that $U$ is an open hereditary set which is not $X$ or the
empty set. Then $K=X\setminus U$ is a nonempty compact subset of
$X$. We claim that \begin{equation}
\bigcup_{i=1}^{N}f_{i}(K)\subseteq K.\label{eq:hered_ifs}\end{equation}
To prove this claim let $x\in K$ and assume that $f_{i}(x)\in U$
for some $i$. Then $x\in r(s^{-1}(f_{i}(x))\subset U$, since $(f_{i}(x),x)$
belongs to the support of $P$. This is a contradiction and the above
inclusion holds. 

We let $F$ be the map defined for all non-empty compact subsets of
$K$ via the formula $F(A)=f_{1}(A)\bigcup\dots\bigcup f_{N}(A)$.
Then we can rewrite equation \eqref{eq:hered_ifs} as $F(K)\subset K$.
It follows that $F^{n}(K)\subseteq K$ for all $n\ge1$. The sequence
$\{F^{n}(K)\}_{n}$ converges to $X$ in the Hausdorff metric (see
\cite{Hutc_80,Bar93,Edg_90}). Thus $X=K$ and this is a contradiction.
\end{proof}
It follows from Theorem \ref{thm:Simplicity} together with Propositions
\ref{pro:IFS_CondL} and \ref{pro:IFS_hereditary} that the $C^{*}$-algebra
$\mathcal{O}(P)$ associated to an iterated function system is always
simple. We record this fact in the following proposition.
\begin{prop}
Suppose that $X$ is a compact nondiscrete metric space and that $(f_{1},f_{2},\dots,f_{N})$
is an iterated function system on $X$. Let $P$ be the Markov operator
associated to this i.f.s. via \eqref{eq:MarkovofIFS}. Then $\mathcal{O}(P)$
is simple.
\end{prop}

\subsection{Collection of continuous maps}

Let $X=[0,1]$, $f_{1}(x)=x$, and $f_{2}(x)=1-x$. We associated
in Example \ref{exa:nonIFS} the following Markov operator \[
P(f)(y)=\frac{1}{2}\sum_{i=1}^{2}f\circ f_{i}(y).\]
Then the support of $P$ is the union of the cographs of $f_{i}$,
$i=1,2$. The formulas defining the topological quiver and the associated
$C^{*}$-correspondence are identical with those for iterated function
systems. The similarities end here, however. Every point $x\in X$
is a base point for at least one loop. All nonreturning loops (see
\cite[Definition 6.5]{MuTo05}) have length at most two. For $x\in X$,
$x\ne1/2$, $(x,x)$ and $((x,1-x),(1-x),x))$ are the only nonreturning
loops. If $x\ne1/2$ the $s^{-1}(x)$ contains exactly two edges,
$(x,x)$ and $(x,1-x)$. If $\alpha$ is a finite path so that $r(\alpha)=1/2$
then $s(\alpha)=1/2$. The topological quiver satisfies condition
(L) since the set of base points of loops with no exits is $\{1/2\}$
and thus it has empty interior in $[0,1]$. The $C^{*}$-algebra $\mathcal{O}(P)$
is not simple, however, because there are many strongly absorbing
closed subsets of $[0,1]$. For example, the sets $\{x,1-x\}$, $x\ne1/2$,
and $\{1/2\}$ are all closed strongly absorbing sets for $P$.

\subsection{Independent Random Variables}

Recall the Markov operator we defined in Example \ref{exa:non_rowfinite}:
$X=[0,1]$, $p(\cdot,y)=m$, Lebesgue measure, for all $x\in[0,1]$.
Then \[
P(f)(y)=\int_{[0,1]}f(x)dm(x),\]
and the support of $P$ is $[0,1]\times[0,1]$. Since $s^{-1}(x)=\{x\}\times[0,1]$
for all $x\in[0,1]$, every point is a base point of a loop and every
loop has an exit. Therefore the topological quiver satisfies condition
(L). Since the only absorbing closed set is $[0,1]$ it follows from
Theorem \ref{thm:Simplicity}  that $\mathcal{O}(P)$ is simple.

\section{From Quivers to Markov Operators\label{sec:From-Quivers}}

We have seen above that a Markov operator determines a topological
quiver. In this section we want to prove that many $C^{*}$-algebras
associated with topological quivers are isomorphic with $C^{*}$-algebras
of specific Markov operators. The key ingredient in our proof is the
so-called \emph{dual topological quiver}.

Let $E=(E^{0},E^{1},r,s,\lambda)$ be a topological quiver, with $E^{0}$
and $E^{1}$ locally compact spaces. We define it's dual as follows.
Set $\widehat{E}^{0}=E^{1}$, $\widehat{E}^{1}=E^{1}*E^{1}$ - the
set of paths of length $2$ in $E$, $\widehat{s}(e_{1},e_{2})=e_{1}$,
$\widehat{r}(e_{1},e_{2})=e_{2}$. We claim that $\widehat{r}$ is
an open map. Let $U_{1}$ and $U_{2}$ be open in $E^{1}$. Since
$\widehat{r}(U_{1}*U_{2})=s^{-1}(r(U_{1}))\bigcap U_{2}$ and $r$
is an open map it follows that $\widehat{r}(U_{1}*U_{2})$ is open
in $\widehat{E}^{0}$. Thus $\widehat{r}$ is open. The family of
Radon measures is defined by $\widehat{\lambda}_{e_{2}}=\lambda_{s(e_{2})}\times\delta_{e_{2}}$,
that is, $\int f(u_{1},u_{2})d\widehat{\lambda}_{e_{2}}(u_{1},u_{2})=\int f(u_{1},e_{2})d\lambda_{s(e_{2})}(u_{1})$.

If $E^{0}$ and $E^{1}$ are compact spaces, $r$ is surjective, and
$\lambda_{v}$ are probability measures for all $v\in E^{0}$, then
we define a Markov operator $P:C(\widehat{E}^{0})\to C(\widehat{E}^{0})$
via\[
Pf(e_{2})=\int f(e)d\lambda_{s(e_{2})}(e).\]
The topological quiver defined by $P$ is $(\widehat{E}^{0},\widehat{E}^{1},\widehat{r},\widehat{s},\widehat{\lambda})$.
We will show that $C^{*}(E)$ and $C^{*}(\hat{E})$ are isomorphic
$C^{*}$-algebra under the assumption that $E$ has no sinks and no
infinite emitters.%
{}
\begin{thm}
\label{thm:Iso-to-dual}Let $\topq$ be a topological quiver with
no sinks and no infinite emitters, that is assume $E_{reg}^{0}=E^{0}$.
If $\widehat{E}=(\widehat{E}^{0},\widehat{E}^{1},\widehat{r},\widehat{s},\lambda)$
is the dual quiver of $E$, then $C^{*}(E)$ and $C^{*}(\widehat{E})$
are isomorphic.\end{thm}
\begin{proof}
Let $(i,\psi)$ and $(\widehat{i},\widehat{\psi})$ be the universal
representations of $C^{*}(E)$ and $C^{*}(\widehat{E})$, respectively.
Let $A=C_{0}(E^{0})$, $\mathcal{X}=\overline{C_{c}(E^{1})}$ be the
$C^{*}$-correspondence associated with $E$, and let $\widehat{A}=C_{0}(\widehat{E}^{0})=C_{0}(E^{1})$
and $\widehat{\mathcal{X}}=\overline{C_{c}(\widehat{E}^{1})}$ be
the $C^{*}$-correspondence associated with $\widehat{E}$. 

An element $f\in C_{b}(E^{1})$ determines an adjointable operator
$T_{f}$ on $\mathcal{X}$ by $T_{f}(\xi)(e)=f(e)\xi(e)$ (see also
\cite[Lemma 3.6]{MuTo05}). Note that $T_{f}^{*}=T_{f^{*}}$. We define
a representation $(\pi,V)$ of $(\widehat{A},\widehat{\mathcal{X}})$
into $C^{*}(E,\lambda)$ by the formulae: $\pi(f)=i^{(1)}(T_{f})$
for $f\in\widehat{A}=C_{0}(E^{1})$ and $V(f_{1}*f_{2})=\psi(f_{1})i^{(1)}(T_{f_{2}})$
for $f_{1},f_{2}\in C_{c}(E^{1})$, where $f_{1}*f_{2}(e_{1},e_{2})=f_{1}(e_{1})f_{2}(e_{2})$.
(Recall that $i^{(1)}$ is the extension of $i$ to $\mathcal{L}(\mathcal{X})$
defined using equation (\ref{eq:PiSuper1}) and the subsequent discussion.)
We claim that $(\pi,V)$ is a Cuntz-Pimsner representation. Let $f\in\widehat{A}$
and $f_{1},f_{2}\in C_{c}(E^{1})$. Then using equation (\ref{eq:Pi1ExtAction})
(with $\pi$ replaced by $i$),\begin{eqnarray*}
V(f\cdot f_{1}*f_{2}) & = & V(T_{f}f_{1}*f_{2})=\psi(T_{f}f_{1})i^{(1)}(T_{f_{2}})\\
 & = & i^{(1)}(T_{f})\psi(f_{1})\widehat{i}^{(1)}(T_{f_{2}})=\pi(f)V(f_{1}*f_{2}),\\
V(f_{1}*f_{2}\cdot f) & = & \psi(f_{1})\widehat{i}^{(1)}(T_{f_{2}f})=\psi(f_{1})i^{(1)}(T_{f_{2}}T_{f})\\
 & = & \psi(f_{1})i^{(1)}(T_{f_{2}})\widehat{i}^{(1)}(T_{f})=V(f_{1}*f_{2})\pi(f).\end{eqnarray*}
If $f_{1},f_{2},g_{1},g_{2}\in C_{c}(E^{1})$ note that $\langle f_{1}*f_{2},g_{1}*g_{2}\rangle_{\widehat{A}}(e)=\overline{f_{2}(e)}\langle f_{1},g_{1}\rangle_{A}(s(e))g_{2}(e)$
and $\langle f_{1},g_{1}\rangle_{A}\circ s\in C_{b}(E^{1})$. Thus
\[
\pi\bigl(\langle f_{1}*f_{2},g_{1}*g_{2}\rangle_{\widehat{A}}\bigr)=\pi(f_{2}^{*})i^{(1)}(T_{\langle f_{1},g_{1}\rangle_{A}\circ s})\pi(g_{2}).\]
Moreover 
\[
i^{(1)}(T_{\langle f_{1},g_{1}\rangle_{A}\circ
  s})\psi(\xi)=\psi(\langle f_{1},g_{1}\rangle_{A}\cdot\xi)=i(\langle
f_{1},g_{1}\rangle_{A})\psi(\xi)=\psi(f_{1})^{*}\psi(g_{1})\psi(\xi).
\]
Thus $i^{(1)}(T_{\langle f_{1},g_{1}\rangle\circ s})=\psi(f_{1})^{*}\psi(g_{1})$.
This implies that\[
V(f_{1}*f_{2})^{*}V(g_{1}*g_{2})=\pi\bigl(\langle f_{1}*f_{2},g_{1}*g_{2}\rangle_{\widehat{A}}\bigr).\]
Thus $(\pi,V)$ is a Toeplitz representation.

We show next that $(\pi,V)$ is a Cuntz-Pimsner covariant representation.
Since $E^{0}=E_{reg}^{0}$ it follows from \cite[Proposition 3.15]{MuTo05}
and \cite[Theorem 3.11]{MuTo05} that if $f\in C_{0}(\widehat{E}_{reg}^{0})$
then $T_{f}\in\mathcal{K}(\mathcal{X})$. Moreover, if $f\in C_{c}(\widehat{E}_{reg}^{0})$
we can use the following construction from the proof of Theorem 3.11
of \cite{MuTo05}. If $K_{f}$ is the support of $f$ there is a finite
cover $\{\, U_{i}\,\}_{i=1}^{n}$ such that $K_{f}\subseteq\bigcup_{i=1}^{n}U_{i}$
and $r|_{U_{i}}:U_{i}\to s(U_{i})$ is a homeomorphism. Moreover we
can choose $U_{i}$ such that $\overline{U}_{i}$ is compact, $s^{-1}(r(\overline{U}_{i}))$
is compact and $r$ restricted to $V_{i}:=s^{-1}(r(U_{i}))$ is a
homeomorphism onto its image (see \cite[Proposition 3.15]{MuTo05}).
If $\{\,\zeta_{i}\,\}_{i=1}^{n}$ is a partition of unity on $K_{f}$
subordinate to $\{\, U_{i}\,\}_{i=1}^{n}$ and $\xi_{i}=f\zeta_{i}^{1/2}$
and $\eta_{i}(\alpha):=\zeta_{i}^{1/2}(\alpha)(\lambda_{r(\alpha)}(\{\alpha\}))^{-1}$
for $\alpha\in E^{1}$ then $\xi_{i},\eta_{i}\in C_{c}(E^{1})$, $\supp\xi_{i}=\supp\eta_{i}\subseteq U_{i}$,
and $T_{f}=\sum_{i=1}^{n}\Theta_{\xi_{i},\eta_{i}}$. Moreover, since
$s^{-1}(r(K_{f}))$ is a closed subset of $\bigcup_{i=1}^{n}s^{-1}(r(\overline{U}_{i}))$
which is compact, it follows that $s^{-1}(r(K_{f}))$ is compact and
$\{\, V_{i}\,\}_{i=1}^{n}$ is an open cover of $s^{-1}(r(K_{f}))$.
Let $\{\,\varsigma_{i}\,\}_{i=1}^{n}$ be a partition of unity subordinate
to $\{\, V_{i}\,\}_{i=1}^{n}$. Then, if $h_{j}=\varsigma_{j}^{1/2}$,
we find that for $f_{1},f_{2}\in C_{c}(E^{1})$ and $(e_{1},e_{2})\in\widehat{E}^{1}$,
\begin{align*}
\widehat{\Phi}(f)(f_{1}*f_{2})(e_{1},e_{2}) & =  f(e_{1})f_{1}(e_{1})f_{2}(e_{2})=\sum_{i,j}f(e_{1})\zeta_{i}(e_{1})f_{1}(e_{1})\varsigma_{j}(e_{2})f_{2}(e_{2})\\
 & =  \sum_{i,j}\xi_{i}(e_{1})\eta_{i}(e_{1})\lambda_{r(e_{1})}(\{e_{1}\})f_{1}(e_{1})h_{j}(e_{2})h_{j}(e_{2})f_{2}(e_{2})\\
 & =  \sum_{i,j}\xi_{i}(e_{1})h_{j}(e_{2})\int_{r^{-1}(r(e_{1}))}\!\!\eta_{_{i}}(u)h_{j}(e_{2})f_{1}(u)f_{2}(e_{2})d\lambda_{s(e_{2})}(u)\\
 & =  \sum_{i,j}\xi_{i}*h_{j}(e_{1},e_{2})\langle\eta_{i}*h_{j},f_{1}*f_{2}\rangle_{\widehat{A}}(e_{2})\\
 & =  \sum_{i,j}\bigl(\xi_{i}*h_{j}\cdot\langle\eta_{i}*h_{j},f_{1}*f_{2}\rangle_{\widehat{A}})(e_{1},e_{2})\\
 & =  \sum_{i,j}\Theta_{\xi_{i}*h_{j},\eta_{i}*h_{j}}(f_{1}*f_{2})(e_{1},e_{2}).\end{align*}
In this computation we used the fact that $r$ is a homeomorphism
on $\supp\xi_{i}=\supp\eta_{i}\subset U_{i}$ and that $r(e_{1})=s(e_{2})$.Thus
\[
\widehat{\Phi}(f)=\sum_{i=1}^{n}\sum_{j=1}^{n}\Theta_{\xi_{i}*h_{j},\eta_{i}*h_{j}}.\]
Therefore\begin{align*}
\pi(f) & =i^{(1)}(T_{f})=\sum_{i=1}^{n}i^{(1)}(\Theta_{\xi_{i},\eta_{i}})=\sum_{i=1}^{n}\psi(\xi_{i})\psi(\eta_{i})^{*}\\
 & =\sum_{i=1}^{n}\sum_{j=1}^{n}\psi(\xi_{i})i^{(1)}(T_{\varsigma_{j}})\psi(\eta_{i})^{*}=\sum_{i,j=1}^{n}\psi(\xi_{i})i^{(1)}(T_{h_{j}})i^{(1)}(T_{h_{j}})^{*}\psi(\eta_{i})^{*}\\
 & =\sum_{i,j=1}^{n}V(\xi_{i}*h_{j})V(\eta_{i}*h_{j})^{*}=\pi^{(1)}\bigl(\widehat{\phi}(f)\bigr).\end{align*}
Using also \cite[Lemma 3.10]{MuTo05} $(\pi,V)$ is a Cuntz-Pimsner
covariant representation. Therefore there exists a $*$-homomorphism
$\pi\rtimes V:C^{*}(\widehat{E})\to C^{*}(E)$ such that $\pi\rtimes V\circ\widehat{i}=\pi$
and $\pi\rtimes V\circ\widehat{\psi}=V$. 

We prove that since $E_{reg}^{0}=E^{0}$, $\pi\rtimes V$ is surjective.
If $a\in C_{0}(E^{0})$ then $a\circ s\in C_{0}(E^{1})$ (\cite[Corollary 3.12]{MuTo05})
and $\phi(a)=T_{a\circ s}$ . Therefore $ $$i(a)=\pi(a\circ s)$
and $i(a)\in\pi(\widehat{A})$. Suppose now that $\xi\in C_{c}(E^{1})$.
Let $\{U_{i}\}_{i=1}^{n}$ be an open cover of the support of $\xi$
and $\{\zeta_{i}\}_{i=1}^{n}$ be a partition of unity subordinate
to $\{U_{i}\}$. Then \begin{eqnarray*}
\psi(\xi) & = & \sum_{i}\psi(\xi)i(\zeta_{i})=\sum_{i}\psi(\xi)i^{(1)}(\phi(\zeta_{i}))\\
 & = & \sum_{i}\psi(\xi)i^{(1)}(T_{\zeta_{i}\circ s})=\sum_{i}V(\xi*\zeta_{i}\circ s).\end{eqnarray*}
Thus $\psi(\xi)\in V(\mathcal{\widehat{X}})$ and it follows that
$\pi\rtimes V$ is surjective.

Finally, we use the gauge invariance uniqueness theorem to show the
injectivity of $\pi\rtimes V$. Let $\gamma_{z}$ and $\beta_{z}$
be the gauge action on $C^{*}(E)$ and $C^{*}(\widehat{E})$, respectively.
Recall that $\gamma_{z}(i(a))=i(a)$ for all $a\in C_{0}(E^{0})$
and $\gamma_{z}(\psi(\xi))=z\psi(\xi)$ for all $\xi\in\mathcal{X}$;
similar statements hold for $\beta_{z}$. Let $f\in C_{0}(E^{1})$.
Then\[
\gamma_{z}\circ\pi\rtimes V(\widehat{i}(f))=\gamma_{z}(\pi(a))=\gamma_{z}(i^{(1)}(T_{f}))\]
which, since $\gamma_{z}(i^{(1)}(T_{f}))(\psi(\xi))=\frac{1}{z}\gamma_{z}\bigl(i^{(1)}(T_{f})\psi(\xi)\bigr)=\frac{1}{z}\gamma_{z}(\psi(T_{f}\xi))=\psi(T_{f}\xi)=i^{(1)}(T_{f})\psi(\xi)$,
equals\[
i^{(1)}(T_{f})=\pi(f)=\pi\rtimes V(\widehat{i}(f))=\pi\rtimes V(\beta_{z}(\widehat{i}(a)).\]
If $f_{1},f_{2}\in C_{c}(E^{1})$\begin{align*}
\gamma_{z}\circ\pi\rtimes V(\widehat{\psi}(f_{1}*f_{2})) & =\gamma_{z}(V(f_{1}*f_{2}))=\gamma_{z}\bigl(\psi(f_{1})i^{(1)}(T_{f_{2}})\bigr)\\
 & =z\psi(f_{1})i^{(1)}(T_{f_{2}})=\pi\rtimes V(z\widehat{\psi}(f_{1}*f_{2}))\\
 & =\pi\rtimes V\circ\beta_{z}(\widehat{\psi}(f_{1}*f_{2}).\end{align*}
Thus $\pi\rtimes V$ is a $*$-isomorphism.\end{proof}
\begin{cor}
\label{cor:Quiver-to-O(P)}Suppose $E=(E^{0},E^{1},r,s,\lambda)$
is a topological quiver with $E^{0}$ and $E^{1}$ compact spaces,
r surjective, and $\lambda_{v}$ a probability measure for all $v\in E^{0}$.
If there are no infinite emitters, then there is a Markov operator
$P$ such that $C^{*}(E)$ is $*$-isomorphic with $\mathcal{O}(P)$.\end{cor}
\begin{rem}
Brenken proved in \cite[Theorem 4.8]{BrB07} that if $E$ is a proper,
range finite topological quiver such that $F_{G}=r(E^{1})-\overline{s(E^{1})}$
closed in $r(E^{1})$, then there is a topological relation $\widehat{E}=(\widehat{E}^{0},\widehat{E}^{1},\widehat{r},\widehat{s},\widehat{\lambda})$
such that $C^{*}(E)$ is isomorphic to $C^{*}(\widehat{E})$. His
topological relation coincides with our dual topological quiver in
the absence of sinks. Our proof is different, though, and our results
imply that the dual topological quiver comes from a Markov operator,
under suitable hypotheses.
\end{rem}
{}
\begin{rem}
Brenken also showed in \cite{BrB07} that without the requirement
that there are no infinite emitters, then Theorem \ref{thm:Iso-to-dual}
can  fail. If $E^{0}=\{v\}$, $E^{1}=[0,1]$, $r(x)=s(x)=v$ for all
$x\in E^{1}$, $\lambda_{v}=m$, Lebesgue measure, then it's dual
topological quiver $\widehat{E}$ is given by the Markov operator
described in Example \ref{exa:non_rowfinite}. Brenken describes these
$C^{*}$-algebras in the comments following Corollary 4.9 of \cite{BrB07}.
He proves that $\Phi(C(E^{0}))\bigcap\mathcal{K}(\mathcal{X})=\emptyset$,
 $C^{*}(E)$ is isomorphic to $\mathcal{O}_{\infty}$, with $K_{0}$
group $\mathbb{Z}$ and trivial $K_{1}$ group, and $C^{*}(\widehat{E})$
is a unital Kirchberg algebra with both $K$ groups $\mathbb{Z}$.
Thus $C^{*}(E)$ and $C^{*}(\widehat{E})$ are not even Morita equivalent,
let alone isomorphic. Nevertheless, one can speculate if Corollary
\ref{cor:Quiver-to-O(P)} is true without the hypothesis that there
are no infinite emitters. The point is that in this case, the Markov
operator will not necessarily come from the dual quiver.
\end{rem}

\bibliographystyle{plain}
\bibliography{MArkovop}

\end{document}